\documentclass[12pt,a4paper]{article}

\usepackage{a4}
\usepackage{url}
\usepackage[margin=3cm]{geometry}

\usepackage{latexsym}

\usepackage[intlimits]{amsmath}
\usepackage{amsfonts}
\usepackage{amssymb}
\usepackage{amsthm}

\usepackage{stmaryrd} 

\usepackage{exscale}

\theoremstyle{plain}

\theoremstyle{definition}

\usepackage[pdftex]{graphicx}
\DeclareGraphicsExtensions{.pdf, .jpg}
\usepackage{color}
\graphicspath{{fig/}}

\definecolor{Blue}{rgb}{0.3,0.3,0.9}
\definecolor{orange}{rgb}{1,0.5,0}
\definecolor{Green}{rgb}{0,0.7,0}

\title{\textbf{Convection-adapted BEM-based FEM}}
\author{Clemens Hofreither\footnote{Johannes Kepler University Linz, Institute for Computational Mathematics, 4040 Linz, Austria} \and Ulrich Langer\footnote{Austrian Academy of Sciences, Johann Radon Institute for Computational and Applied Mathematics, 4040 Linz, Austria} \and Steffen Wei\ss er\footnote{Saarland University, Department of Mathematics, 66041 Saarbr\"ucken, Germany}}

\usepackage[final,pdftex,
      pdftitle={Convection-adapted BEM-based FEM},
      pdfauthor={Clemens Hofreither and Ulrich Langer and Steffen Weisser},
      pdfkeywords={Convection-diffusion-reaction problems; non-standard finite element methods; BEM-based FEM; Local Trefftz methods},
	colorlinks, 
	bookmarksnumbered, 
	bookmarksopen, 
	pdfstartview=FitH, 
	linkcolor=black, 
	citecolor=black, 
	urlcolor=black, 
	filecolor=black %
]{hyperref}
\newcommand{\norm}[1]{\|#1\|}
\newcommand{\sprod}[2]{\langle #1,\, #2 \rangle}

\begin{document}

\maketitle

\begin{abstract}
We present a new discretization method for homogeneous convection-diffusion-reaction 
boundary value problems in 3D that is a non-standard finite element method 
with PDE-harmonic shape functions on polyhedral elements.
The element stiffness matrices are constructed by means of 
local boundary element techniques.
Our method, which we refer to as a BEM-based FEM, can therefore be considered
a local Trefftz method with element-wise (locally) PDE-harmonic shape functions.
The Dirichlet boundary data for these shape functions
is chosen according to a convection-adapted procedure
which solves projections of the PDE onto the edges and faces of the elements.
This improves the stability of the discretization method for convection-dominated problems
both when compared to
a standard FEM and to previous BEM-based FEM approaches,
as we demonstrate in several numerical experiments.
\\[0.5ex]
{\footnotesize
\noindent\textbf{Keywords} Convection-diffusion-reaction problems $\cdot$ 
  non-standard finite element methods $\cdot$
  BEM-based FEM $\cdot$
  Local Trefftz methods
\\[0.5ex]
\noindent\textbf{Mathematics Subject Classification (2000)} 
65N30, 65N38
}
\end{abstract}

\section{Introduction}
\label{sec:Introduction}
The BEM-based FEM was introduced in \cite{CopelandLangerPusch2009} 
on the basis of ideas borrowed from boundary element 
domain decomposition methods originally proposed 
by G.~C.~Hsiao and W.~L.~Wendland in \cite{HsiaoWendland:1991a}.
This class of discretization methods uses PDE-harmonic shape functions
in every element of a polyhedral mesh. In order to generate the local 
stiffness matrices efficiently, boundary element techniques are 
employed locally. This is the reason why these non-standard Finite Element Methods 
are called BEM-based FEM. A BEM-based FEM can also be considered a local  Trefftz FEM.
The papers  \cite{HofreitherLangerPechstein2010} and   \cite{Hofreither2011} provide the a priori discretization error analysis  
with respect to the energy and $L_2$ norms, respectively,
where homogeneous diffusion problems serve as model problems.
In \cite{RjasanowWeisser2014}, a new construction of the approximation space was proposed, which employs the polygonal faces of the polyhedral elements.
Residual-type a~posteriori discretization error estimates were derived in
\cite{Weisser2011} and extended in \cite{Weisser2015Enumath}. These a~posteriori discretization error estimates
can be used to derive adaptive versions of the BEM-based FEM,
see also the PhD thesis by S.~Wei{\ss}er \cite{WeisserPHD}.
In addition to the low order approximation techniques, high order trial functions were introduced, discussed and studied in \cite{RjasanowWeisser2012,Weisser2014,Weisser2014PAMM}, which open the development towards fully $hp$-adaptive strategies.
Fast FETI-type solvers for solving the large linear systems arising 
from the BEM-based FEM discretization of diffusion problems 
were studied in \cite{HofreitherLangerPechstein:2014a}.
Furthermore, the ideas of the BEM-based FEM are transferred into other application areas. There are, for example, first results on vector valued, $H(\operatorname{div})$-conforming approximations \cite{EfendievGalvisLazarovWeisser2014} and on time dependent problems \cite{Weisser2014WCCM}.

The use of PDE-harmonic shape functions seems to 
be especially appropriate for convection-diffusion 
problems. The first results in this direction were presented 
in \cite{HofreitherLangerPechstein2011}, see also the 
PhD thesis by C. Hofreither \cite{HofreitherPHD}.
The shape functions used in these works are PDE-harmonic 
in the interior of the polyhedral elements, but their traces 
on the boundaries of the polyhedral elements are still 
piecewise linear and not adapted to the convection.
There is a close relation between this BEM-based FEM with piecewise
linear boundary data and the so-called method of residual-free bubbles
\cite{BrezziEtAl:1992,BrezziEtAl:1997,BrezziEtAl:1998,FrancaEtAl:1998,BrezziEtAl:1999}.
Indeed, it has been shown in \cite{HofreitherPHD} that the BEM-based FEM,
with exact evaluation of the Steklov-Poincar\'e operator, is equivalent
to the method of residual-free bubbles with exactly computed bubbles.
Since the latter has been shown to be a stable method for convection-dominated problems,
it seems clear that also the BEM-based FEM should have advantageous stability properties.
It should be noted that neither the Steklov-Poincar\'e operator nor the computation of the
residual-free bubbles can be realized exactly in practice.
However, the early numerical experiments in \cite{HofreitherLangerPechstein2011,HofreitherPHD}
demonstrate the stabilizing properties of the BEM-based FEM.

In this paper, we aim to further enhance this stabilizing effect.
Our approach is to construct a new BEM-based FEM for 
convection-diffusion-reaction boundary value problems
which employs
basis functions that are still PDE-harmonic within the elements,
and at the same time convection-adapted on the element boundaries.
This adaption to the convection is obtained by solving 
projected 1D and 2D boundary value problems at all edges and 
faces of the polyhedral elements, respectively.
The projected 1D convection-diffusion-reaction problems on the edges 
can be solved analytically, whereas the 2D face problems are solved numerically
by means of a Streamline Upwind/Petrov-Galerkin (SUPG) finite element method 
on an auxiliary face triangulation.
This approach extends the stable applicability of 
the new method considerably, as is shown by our 
numerical experiments.

The remainder of this paper is structured as follows. 
In Section~2, we derive the skeletal
variational formulation that will be the starting point for the discretization.
Section~3 provides the boundary integral representations of local Steklov-Poincar\'{e} operators.
The main results are contained in Section~4, where we construct the convection-adapted 
PDE-harmonic shape functions
and derive a fully discretized skeletal variational formulation as
linear systems of algebraic equations.
In Section~5, we present and discuss some numerical results illustrating 
that the convection-adapted BEM-based FEM works well in the convection-dominated case. 
Finally, Section~6 draws some conclusions and provides an outlook on further work.

\section{Derivation}

In this section, we briefly derive a so-called skeletal variational formulation
for a convection-diffusion-reaction problem.
As our model problem, we consider the pure Dirichlet boundary value problem
\begin{equation}
    \begin{aligned}
      L u = -\operatorname{div}(A \nabla u) + b \cdot \nabla u + cu &= 0     \quad\text{in } \Omega, \\
      u & = g   \quad\text{on } \partial\Omega
    \end{aligned}
    \label{eq:pde}
\end{equation}
in a bounded Lipschitz domain $\Omega \subset \mathbb R^3$.
Here $A(x) \in \mathbb R^{3 \times 3},$ $b(x) \in \mathbb R^3$, and $c(x) \in \mathbb R$ are the coefficient
functions of the partial differential operator $L$,
and $g \in H^{1/2}(\partial\Omega)$ is the given Dirichlet data.
We assume that $A(\cdot)$ is symmetric and uniformly positive, 
and that $c(\cdot)$ is non-negative.
The corresponding variational formulation reads as follows:
find $u \in H^1(\Omega)$ with $\gamma^0_\Omega u = g$ such that
\begin{equation}
  \label{eq:standardvf}
  \int_\Omega \left(A \nabla u\cdot\nabla v + b \cdot \nabla u \,v + cuv \right) dx = 0
  \qquad \forall\, v \in H^1_0(\Omega),
\end{equation}
where $\gamma^0_\Omega : H^1(\Omega) \rightarrow H^{1/2}(\partial\Omega)$
refers to the Dirichlet trace operator from the domain $\Omega$ to its boundary,
and $H^1_0(\Omega) = \{ v \in H^1(\Omega) : \gamma^0_\Omega v = 0 \}$.
We require that the coefficients $A$, $b$, $c$ are $L^\infty(\Omega)$
and that there exists a unique solution of \eqref{eq:standardvf}.

We assume that we have a finite decomposition $\mathcal T$ of $\Omega$ into mutually disjoint Lipschitz polyhedra.
These will play the role of elements as in a finite element method, 
but we do not require the existence of a reference element to which the elements can be mapped
and instead allow $\mathcal T$ to contain an arbitrary mixture of polyhedral element shapes.
We require that the coefficients 
$A(\cdot),$ $b(\cdot),$ and $c(\cdot)$
are piecewise
constant with respect to the polyhedral mesh $\mathcal T$.

It follows from the variational formulation and the density of $C^\infty_0(T)$ in $L_2(T)$ that
$\operatorname{div}(A \nabla u) = b \cdot \nabla u + cu \in L_2(T)$ for every element $T \in \mathcal T$.
Therefore, the flux $A \nabla u$ is in $H(\operatorname{div}, T)$.
Let $n_T$ denote the outward unit normal vector on $\partial T$.
Then the flux has a well-defined normal trace $\gamma_T^1 u = A\nabla u \cdot n_T \in H^{-1/2}(\partial T)$,
also called the \emph{conormal derivative} of $u$ (cf.~\cite{GiraultRaviart:1986a}).
Moreover, we have the generalized Green's identity
\begin{align}
\label{eq:Green}
  \int_T A\nabla u \cdot \nabla v\, dx
  = - \int_T \operatorname{div}(A\nabla u) v\, dx + \langle \gamma^1_T u ,\, \gamma^0_T v \rangle
  \qquad \forall v \in H^1(T),
\end{align}
where $\langle\cdot,\cdot\rangle$ denotes the duality pairing on 
$H^{-1/2}(\partial T) \times H^{1/2}(\partial T)$.
We mention that the particular element boundary $\partial T$ will always be clear by context.
Inserting \eqref{eq:Green} into \eqref{eq:standardvf} and 
recalling
that $Lu = 0$ in $L_2(T)$,
we obtain
\begin{align*}
  0
  = \sum_{T\in \mathcal T} \int_T \left(A \nabla u\cdot\nabla v + b \cdot \nabla u \,v + cuv \right) dx
  = \sum_{T\in \mathcal T} \Big( \underbrace{\int_T Lu\,v \,dx}_{=0}
    + \langle \gamma^1_T u ,\, \gamma^0_T v \rangle \Big).
\end{align*}
Fix now some element $T \in \mathcal T$ and observe that $u\arrowvert_T$ is the unique solution of the
local boundary value problem
\begin{equation}
    \label{eq:localbvp}
   \text{find } \varphi \in H^1(T): \quad
   L \varphi = 0, \qquad \gamma^0_T \varphi = \gamma^0_T u.
\end{equation}
We require that these local problems have unique solutions too.
Denoting by $S_T : H^{1/2}(\partial T) \rightarrow H^{-1/2}(\partial T)$ the
\emph{Steklov-Poincar\'e operator} or \emph{Dirichlet-to-Neumann map} for this local problem,
we therefore have $\gamma^1_T u = S_T \gamma^0_T u$.
Writing $u_{\partial T} = \gamma^0_T u$ and analogously $v_{\partial T}$ for the traces 
onto the element boundary,
we get
\[
   \sum_{T \in \mathcal T} \langle S_T u_{\partial T},\, v_{\partial T} \rangle = 0
   \qquad \forall v \in H^1_0(\Omega).
\]

Note that all terms in the last formulation are defined only on the element boundaries $\partial T$.
Let $\Gamma_S = \bigcup_T \partial T$ denote the \emph{skeleton} of the mesh $\mathcal T$,
 and let the
skeletal function space $W = H^{1/2}(\Gamma_S)$ consist of the traces of all functions from $H^1(\Omega)$
on $\Gamma_S$.
Then we are looking for a skeletal function $u \in W$
which satisfies the Dirichlet boundary condition $u\arrowvert_{\partial\Omega}=g$
and the skeletal variational formulation
\begin{equation}
  \label{eq:skelvf}
  \sum_{T \in \mathcal T} \langle S_T u_{\partial T} ,\, v_{\partial T} \rangle = 0
   \qquad \forall\, v \in W_0 = \left\{v \in W: v\arrowvert_{\partial\Omega} = 0\right\}.
\end{equation}
This skeletal variational formulation is equivalent to the
standard variational formulation \eqref{eq:standardvf}
in the sense that the traces
$u_{\partial T} \in H^{1/2}(\partial T)$ obtained from \eqref{eq:skelvf} match the traces
$\gamma^0_T u$ of the function $u \in H^1(\Omega)$ obtained from \eqref{eq:standardvf}.
Conversely, $u\arrowvert_T$ can be recovered from $u_{\partial T}$
by solving a local Dirichlet problem of type \eqref{eq:localbvp} in $T$. 
We will call the solution of this local problem the PDE-harmonic extension of $u_{\partial T}$.
We remark that another interpretation of \eqref{eq:skelvf} is that of a weak enforcement
of the continuity of conormal derivatives on inter-element boundaries.

\section{Boundary integral operators}

Evaluating the Dirichlet-to-Neumann map $S_T$ used above essentially corresponds to solving a local problem $L\varphi=0$ on $T$
with the given Dirichlet data and then obtaining the conormal derivative $\gamma^1_T \varphi$ of its solution.
Even these local problems are in general not analytically solvable and require numerical approximation.
At the core of our method, we approximate $S_T$ using a boundary element 
technique
obtained by the Galerkin discretization of element-local boundary integral equations.
Some standard results on boundary integral equations are outlined in the following.
A more detailed treatment of these topics can be found in, e.g., 
\cite{HsiaoWendland:2008a,Mclean2000,SauterSchwab:eng,Steinbach:numell:eng}.

An important prerequisite for defining the boundary integral operators is the knowledge of
a \emph{fundamental solution}.
A fundamental solution of the partial differential operator $L$
is a function $G(x,y)$ such that $L_x G(x,y) = \delta(y-x)$,
where $\delta$ is the Dirac $\delta$-distribution and $x$, $y \in \mathbb R^d$.
A fundamental solution for $L$ from \eqref{eq:pde}
with constant coefficients $A$, $b$, $c$
is given in \cite{SauterSchwab:eng}.
In fact,
in $\mathbb R^3$ and under the assumption $c + \norm{b}_{A^{-1}}^2 \ge 0$, we have
\[
  G(x,y) = \frac 1 {4\pi \sqrt{\det A}} \frac{\exp \left( b^\top A^{-1} (x-y)
      - \lambda \norm{x-y}_{A^{-1}} \right)}{\norm{x-y}_{A^{-1}}}\,,
\]
where $\norm{x}_{A^{-1}} = \sqrt{x^\top A^{-1} x}$ and $\lambda = \sqrt{c + \norm{b}_{A^{-1}}^2}$.

For our setting, we more generally assume that the coefficients $A$, $b$, $c$
are constant only within each element.
This leads to a potentially different fundamental solution
in each element $T$, in the following denoted by $G_T(x,y)$,
and allows us to treat PDEs with piecewise constant coefficients.

We now introduce the local boundary integral operators
\begin{gather*}
  V_T : H^{-1/2}(\partial T) \rightarrow H^{1/2}(\partial T), \quad
  K_T : H^{1/2}(\partial T) \rightarrow H^{1/2}(\partial T),   \\
  K'_T: H^{-1/2}(\partial T) \rightarrow H^{-1/2}(\partial T), \quad
  D_T : H^{1/2}(\partial T)  \rightarrow H^{-1/2}(\partial T),
\end{gather*}
called, in turn, the \emph{single layer potential}, \emph{double layer potential},
\emph{adjoint double layer potential}, and \emph{hypersingular} operators.
For sufficiently regular arguments, they admit the integral representations
\begin{align*}
  (V_Tv)(y) &= \int_{\partial T} G_T(x,y) v(x) \, ds_x,             \\
  (K_Tu)(y) &= \int_{\partial T} \widetilde{\gamma^1_{T,x}} G_T(x,y) u(x) \, ds_x,    \\
  (K'_Tv)(y) &= \int_{\partial T} \gamma^1_{T,y} G_T(x,y) v(x) \, ds_x,   \\
  (D_Tu)(y) &= -\gamma^1_{T,y} \int_{\partial T} \widetilde{\gamma^1_{T,x}} G_T(x,y) \big(u(x) - u(y)\big) \, ds_x,
\end{align*}
where $\gamma^1_{T,y}$ refers to the conormal derivative $\gamma^1_T$ with respect to the variable $y$,
whereas
$\widetilde{\gamma^1_{T,x}}$ refers to the modified conormal derivative with respect to 
$x$, i.e.,
\[
 \widetilde{\gamma^1_T} u = \gamma^1_T u  +  (b \cdot n_T) \gamma^0_T u,
\]
which is associated with the adjoint problem.
We have the following two representations
of the Steklov-Poincar\'e operator in terms of the boundary integral operators:
\begin{equation}
  \label{eq:steklov}
  S_T = V_T^{-1} (\tfrac12 I + K_T)
      = D_T + (\tfrac12 I + K_T') V_T^{-1} (\tfrac12 I + K_T).
\end{equation}

\section{Discretization}
\label{sec:Discretization}

\subsection{Discretization of the skeletal function space}
\label{subsec:DiscretizationSkeletalFunctionSpace}

We employ a Galerkin approach to the discretization of the skeletal variational 
formulation \eqref{eq:skelvf}. 
To this end, we notice that every element boundary $\partial T$ is composed of 
open polygonal faces $\mathcal F_T$, straight edges $\mathcal E_T$ and 
nodes $\mathcal N_T$ located in the corner points of the element. 
Neighboring elements either share a common face, edge or node. 
This gives a natural description of the skeleton 
${\overline {\mathcal F}}=\bigcup_{T\in\mathcal T}{\overline {\mathcal F}}_T$ 
of the decomposition~$\mathcal T$.
In the following, we construct a discrete trial space $W_h\subset W$ over $\mathcal F$. 
The PDE-harmonic extensions of its basis functions can thus be interpreted as three 
dimensional trial functions for the approximation of $u$.

If the skeleton~$\mathcal F$ consists only of triangular faces or if it is triangulated, 
a straightforward choice for $W_h$ is the space of piecewise linear and globally continuous 
functions on~$\mathcal F$. This strategy was introduced in~\cite{HofreitherLangerPechstein2011} 
for the convection-diffusion-reaction equation. 
An enhanced strategy was proposed in~\cite{RjasanowWeisser2014} for the pure diffusion equation. 
It makes use of a hierarchical construction of trial functions on polyhedral elements with polygonal faces.
Extending this idea, we define the basis functions of $W_h$ with the help of 
PDE-harmonic extensions on edges and faces, which might be triangles or of polygonal shape.

More precisely, 
for each node $z_i\in\mathcal N=\bigcup_{T\in\mathcal T}\mathcal N_T$,
we introduce a skeletal basis function $\varphi_i\in W$ satisfying
\begin{gather}
 \varphi_i(z_j) = \delta_{ij}\quad\mbox{for } z_j\in\mathcal N,\nonumber\\
 L_E\varphi_i = 0 \quad\mbox{on } E\in\mathcal E,\label{eq:defBasisFkt}\\
 L_F\varphi_i = 0 \quad\mbox{on } F\in\mathcal F,\nonumber
\end{gather}
where $L_E$ and $L_F$ are projections of the differential operator~$L$ onto the edge~$E$ and face~$F$, 
respectively, and $\mathcal E=\bigcup_{T\in\mathcal T}\mathcal E_T$. 
Thus, the functions $\varphi_i$ are defined implicitly as local solutions 
of boundary value problems on edges and faces of the decomposition. 
Equivalently, one can say that these functions are defined via PDE-harmonic extensions. The nodal data is first extended $L_E$-harmonically along the edges and afterwards, the data on the edges is extended into the faces with the help of a $L_F$-harmonic operator. 
 
For the definition of $L_E$ and $L_F$, let $F\in\mathcal F$ be a face and $E\in\mathcal E$ an edge on the boundary of $F$. By rotation and translation of the coordinate system, we map the face $F$ into the $(e_1,e_2)$-plane and the edge $E$ onto the $e_1$-axis of the Euclidean coordinate system $(e_1,e_2,e_3)$ such that one node of $E$ lies in the origin.
Thus, we have an orthogonal matrix $B \in \mathbb R^{3\times3}$ and a vector $d \in \mathbb R^3$
such that
\[
    \tilde x \mapsto x = B\tilde x+d\quad\mbox{and}\quad \tilde\varphi(\tilde x) = \varphi(B\tilde x+d),
\]
and the differential equation in~\eqref{eq:pde} yields
\begin{equation}\label{eq:MappedPDE}
  -\operatorname{div}(A \nabla \varphi) + b \cdot \nabla \varphi + c\varphi 
  = -\operatorname{div}_{\tilde x}(BAB^\top \nabla_{\tilde x} \tilde\varphi) + Bb \cdot \nabla_{\tilde x} \tilde\varphi + c\tilde\varphi 
  = 0.
\end{equation}
Furthermore, we assume for the definition of trial functions that they only vary in tangential direction to the face and edge, respectively, i.e.,
$\frac{\partial\tilde\varphi}{\partial\tilde x_3}=0$ in $F$ and
$\frac{\partial\tilde\varphi}{\partial\tilde x_2}=\frac{\partial\tilde\varphi}{\partial\tilde x_3}=0$ on $E$.
Therefore, the dependence in~\eqref{eq:MappedPDE} reduces to two and one coordinate directions such that $L_F$ and $L_E$ are defined as differential operators in two and one dimensions using the described coordinate system. 
Overall, the basis functions are constructed with the help of the convection-diffusion-reaction equation on the faces and edges, where the diffusion matrix and the convection vector are adjusted in a proper way.

Attentive readers have noticed that $A$, $b$ and $c$ are assumed to be constant on each polygonal subdomain such that they are not well-defined on faces and edges. Therefore,
in order to make the operators $L_E$ and $L_F$ well-defined,
we have to understand $A$, $b$ and $c$ in~\eqref{eq:MappedPDE} on $F$ and $E$ as averaged quantities over the neighboring polyhedra of the face~$F$ and the edge~$E$, respectively. To simplify notation, we will omit the coordinate transformation in the following and abbreviate the transformed diffusion matrix~$BAB^\top$ and convection vector~$Bb$ to $A_F$ and $b_F$, respectively. Furthermore, we will treat the basis functions~$\varphi_i$ as functions of two or one variable depending on the underlying domain $F$ or $E$.

Having defined the basis functions $\varphi_i$, we obtain a discretization space
$W_h = \operatorname{span} \{ \varphi_i :i=1,\ldots,|\mathcal N|\} \subset W$,
where $|\mathcal N|$ denotes the number of nodes in the polyhedral mesh.
Assuming that the given Dirichlet data $g$ can be extended to the skeleton by a function in $W_h$,
we thus arrive at the following Galerkin equations as the discrete version of \eqref{eq:skelvf}:
find $u_h \in W_h$ such that $u_h\arrowvert_{\partial\Omega} = g$ and
\begin{equation}
  \sum_{T \in \mathcal T} \sprod {S_T u_{h,\partial T}} {v_{h,\partial T}} = 0
  \qquad \forall \, v_h \in W_{h,0} = W_h \cap W_0.
  \label{eq:discrvf}
\end{equation}
In the general case, the Dirichlet data $g$ can be approximated by the introduced trial functions $\varphi\in W_h$ using interpolation (if continuous) or $L_2$-projection.

\subsection{Approximation of the skeletal basis functions}
\label{subsec:ApproximationSkeletalFunctionSpace}

The construction of the skeletal basis functions $\{ \varphi_i \}$ involves the solution of
certain lower-dimensional boundary value problems which cannot, in general, be done exactly.
In this section, we therefore construct computable approximations to the exact basis functions.

In~\eqref{eq:defBasisFkt}, $L_E$ describes an ordinary differential operator of second order with constant and scalar-valued coefficients. Thus the boundary value problems on the edges can be solved analytically and the restrictions of the functions in $W_h$ to each edge $E\in\mathcal E$ can be written in closed form.

\begin{figure}[tb]
 \centering
 \includegraphics[trim=4cm 2cm 4cm 4cm, width=0.3\textwidth, clip]{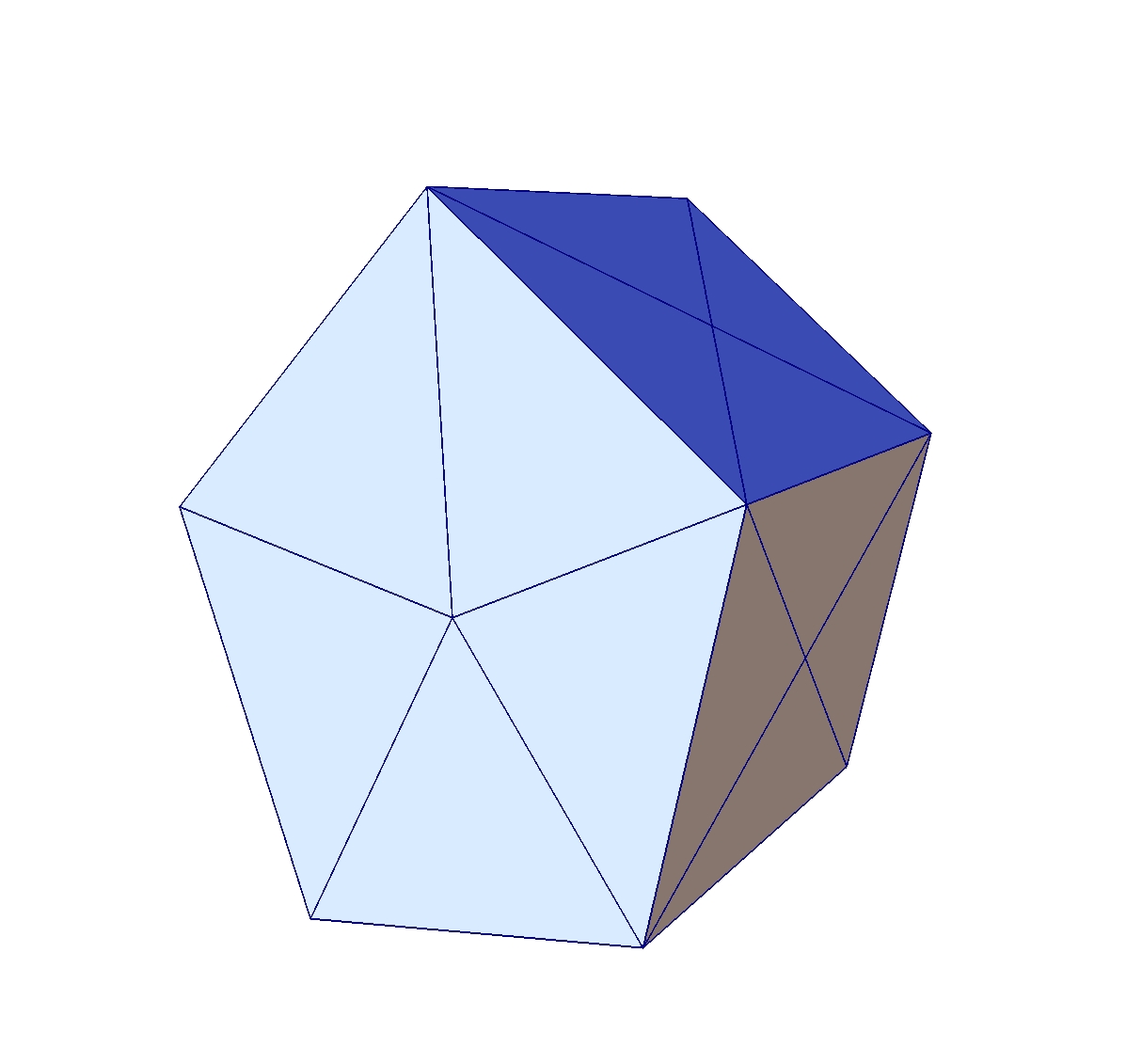}\hfil
 \includegraphics[trim=4cm 3cm 4cm 4cm, width=0.3\textwidth, clip]{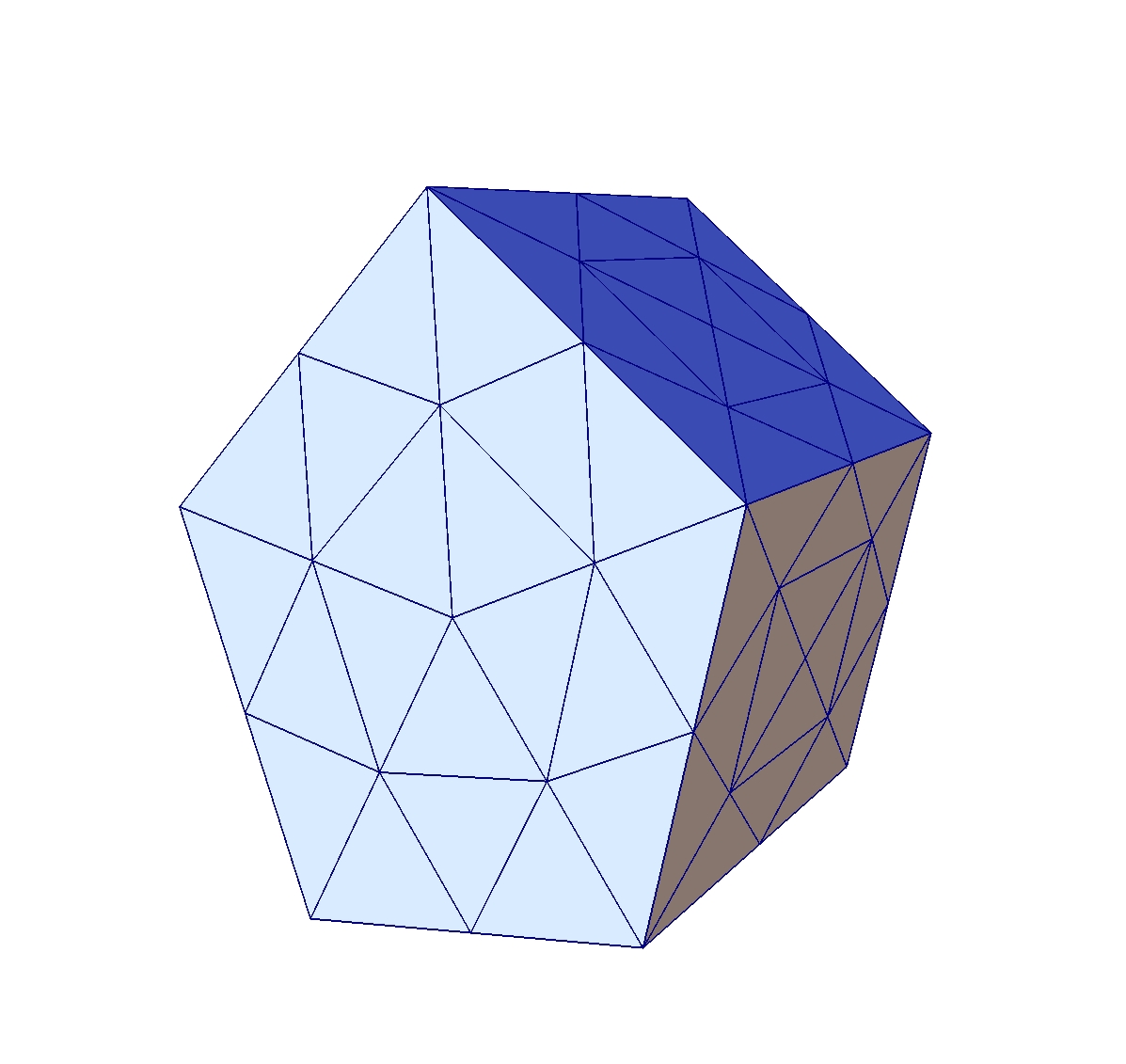}\hfil
 \includegraphics[trim=4cm 3cm 4cm 4cm, width=0.3\textwidth, clip]{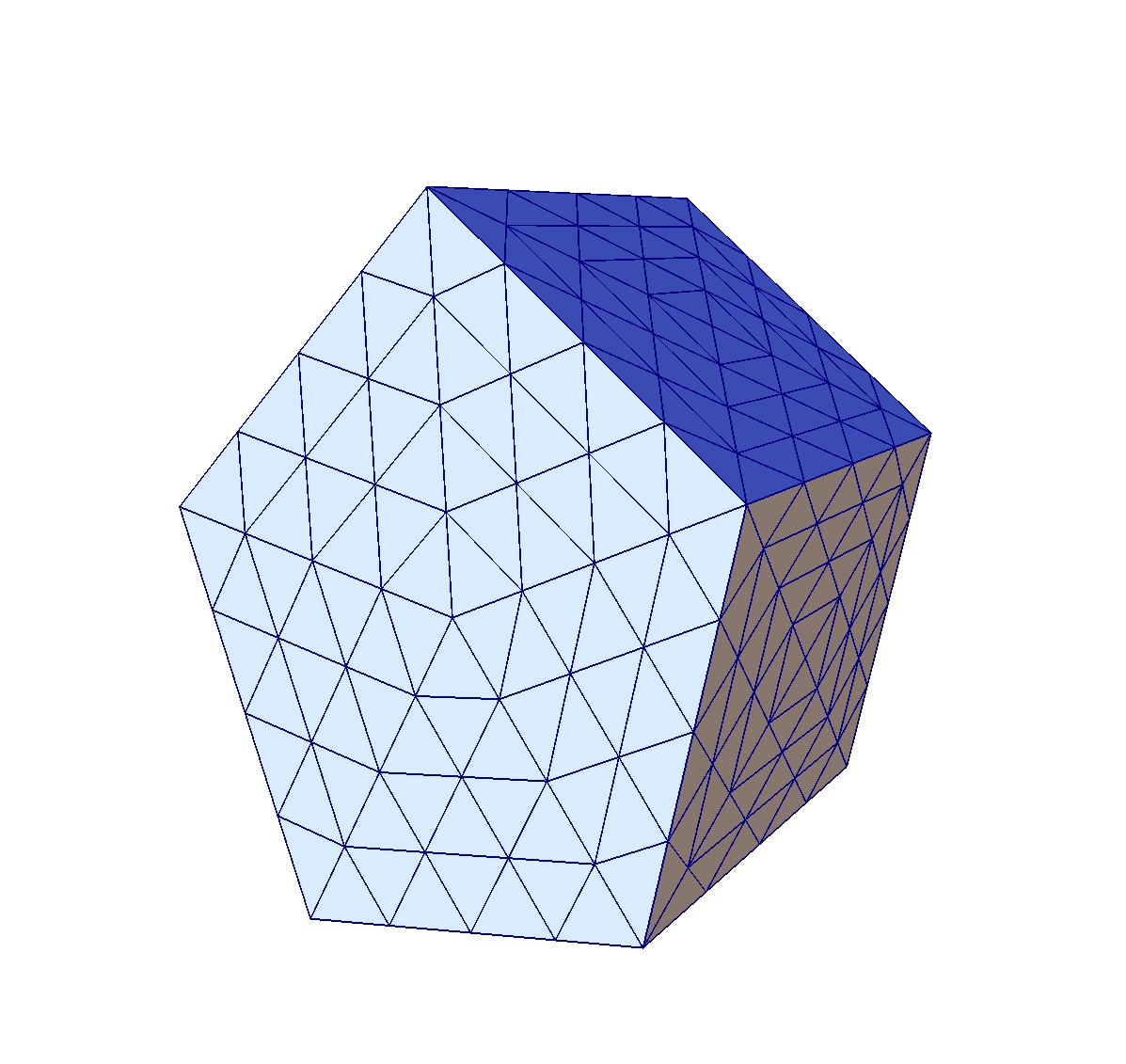}
 \caption{Auxiliary triangulation $\mathfrak T_\ell(\mathcal F_T)$ with $\ell=0,1,2$ for the surface of $T\in\mathcal T$.}
 \label{fig:AuxiliaryTriangulation}
\end{figure}

The two-dimensional problems involving the operator $L_F$ on the faces, however, need to be solved approximately.
We choose a SUPG method, see~\cite{BrooksHughes1982}, on each face since we might have convection-dominated problems there. To this end, we introduce an auxiliary triangulation of each face in such a way that the meshes are matching on common edges as well as quasi-uniform and shape-regular in the usual sense, with constants which are uniform over all faces $F\in\mathcal F$. 
In the case of elements and faces which are star-shaped with respect to a ball and a circle, respectively, we can use, for example, the construction in~\cite{RjasanowWeisser2014}. Here, a first coarse triangulation of a face~$F\in\mathcal F$ is obtained by connecting its nodes with the center of the inscribed circle. This auxiliary mesh is denoted by $\mathfrak T_0(F)$. Afterwards, the meshes $\mathfrak T_\ell(F)$ of level $\ell \geq 1$ are defined recursively by splitting each triangle of the previous level into four similar triangles
by connecting its edge midpoints.
This strategy yields a conforming triangulation $\mathfrak T_\ell(\mathcal F_T)=\bigcup_{F\in\mathcal F_T}\mathfrak T_\ell(F)$ of the surface of the polyhedral elements, see Figure~\ref{fig:AuxiliaryTriangulation}, as well as a conforming triangulation $\mathfrak T_\ell(\mathcal F)=\bigcup_{F\in\mathcal F}\mathfrak T_\ell(F)$ of the whole skeleton.
Furthermore, it also induces a discretization of each edge $E \in \mathcal E$ into line segments of equal size
which we denote by $\mathfrak T_\ell(E)$.

Let $\mathfrak W_\ell(\mathcal F)$ denote the space of piecewise linear functions over the auxiliary triangulation $\mathfrak T_\ell(\mathcal F)$ which are continuous on the skeleton. We denote the restrictions of $\mathfrak W_\ell(\mathcal F)$ onto a face $F$ and the boundary $\partial T$ of an element by $\mathfrak W_\ell(F)$ and $\mathfrak W_\ell(\mathcal F_T)$, respectively. Furthermore, let $\mathfrak W_{\ell,0}(F)$ be the subspace of $\mathfrak W_\ell(F)$ which contains only functions vanishing on the boundary of~$F$.

Now, we have the ingredients to handle the approximation of the basis functions on the faces.
For any $\varphi_i\in W_h$, we define its approximation $\varphi^i_\ell\in\mathfrak W_\ell(F)$ on each face $F\in\mathcal F$ with $z_i\in\partial F$
such that it coincides with the piecewise linear (with respect to $\mathfrak T_\ell(E)$) nodal interpolation of $\varphi_i$ on the edges of $F$
and such that it fulfills the SUPG formulation: find $\varphi^i_\ell\in\mathfrak W_\ell(F)$ such that boundary conditions on $\partial F$, as described above, are fulfilled and
\begin{multline}\label{eq:VFonFace}
  \int_F(A_F\nabla\varphi^i_\ell\cdot\nabla\phi
         + b_F\cdot\nabla\varphi^i_\ell \,\phi
         + c_F\varphi^i_\ell \,\phi)\\
    + \delta_F\int_F(b_F\cdot\nabla\varphi^i_\ell \, b_F\cdot\nabla\phi
               + c_F\varphi^i_\ell \, b_F\cdot\nabla\phi)
  = 0 \quad\forall \phi\in\mathfrak W_{\ell,0}(F),
\end{multline}
where $\delta_F\geq0$ is a stabilization parameter which is set to zero in the diffusion-dominated case. On all faces $F\in\mathcal F$ with $z_i\notin\partial F$, it is $\varphi_\ell^i \equiv 0$.
Finally, we obtain the approximate trial space as $W_\ell=\operatorname{span}\{ \varphi^i_\ell:i=1,\ldots,|\mathcal N|\}
\subset \mathfrak W_\ell(\mathcal F) \subset W$.

Assuming that the given Dirichlet data $g$ can be extended to the skeleton by a function in $W_\ell$,
we thus arrive at the following Galerkin equations as the approximated version of \eqref{eq:discrvf}:
find $u_\ell \in W_\ell$ such that $u_\ell\arrowvert_{\partial\Omega} = g$ and
\begin{equation}
  \sum_{T \in \mathcal T} \sprod {S_T u_{\ell,\partial T}} {v_{\ell,\partial T}} = 0
  \qquad \forall \, v_\ell \in W_{\ell,0} = W_\ell \cap W_0.
  \label{eq:approxdiscrvf}
\end{equation}
In the general case, the Dirichlet data $g$ can be approximated by the introduced trial functions $\varphi_\ell^i\in W_\ell$ on the boundary faces $F\subset\partial\Omega$ using
interpolation (if continuous) or $L_2$-projection.

\subsection{Discretization of the Dirichlet-to-Neumann map}
\label{subsec:DiscD2N}

The Dirichlet-to-Neumann maps $S_T$ in \eqref{eq:approxdiscrvf} will be evaluated via the
boundary integral operator representation \eqref{eq:steklov}. 
In the case of diffusion-reaction problems the symmetric representation is preferred, since it preserves the symmetry.
In this paper, however, we use the first representation since we deal with a non-symmetric problem. Furthermore, this representation is easier and can be implemented more efficiently.
This formula still contains the inverse of the single layer potential operator $V_T$, which is in general not computable exactly.
Hence, we also need to approximate the bilinear form $\sprod {S_T \cdot}{\cdot}$.
To do this, we employ a mixed continuous
piecewise linear/piecewise constant scheme, where Dirichlet data are approximated linearly,
while Neumann data are approximated by piecewise constant functions,
as described in, e.g., \cite{Costabel:1987a, Steinbach:numell:eng, HsiaoSteinbachWendland:ddbem}.

Let $\phi_{T,i}\in\mathfrak W_\ell(\mathcal F_T)$ denote the nodal piecewise linear functions
restricted to the local mesh $\mathfrak T_\ell(\mathcal F_T)$,
where now $i$ enumerates the vertices of the auxiliary triangulation $\mathfrak T_\ell(\mathcal F_T)$.
Furthermore, introduce a space of piecewise (per triangle of $\mathfrak T_\ell(\mathcal F_T)$) constant boundary functions
spanned by the basis $\{\psi_{T,k}\}$, where $k$ enumerates the triangles $\tau_k \in \mathfrak T_\ell(\mathcal F_T)$,
such that $\psi_{T,k} \equiv 1$ on $\tau_k$ and $\psi_{T,k} \equiv 0$ on all other triangles.

For any function $u_T \in H^{1/2}(\partial T)$, its corresponding Neumann data
can be written, according to \eqref{eq:steklov}, as
\[
  t_T = S_T u_T = V_T^{-1} (\tfrac12 I + K_T) u_T \in H^{-1/2}(\partial T).
\]
For any piecewise linear function $u_{T\ell} \in \mathfrak W_\ell(\mathcal F_T)=\operatorname{span}\{ \phi_{T,i} \}$, we now compute an approximation $t_{T\ell} \approx t_T$
of its Neumann data by the Galerkin projection of the equation
$V_T t_T = (\tfrac12 I + K_T) u_T$ to the piecewise constant functions.
In other words, we seek $t_{T\ell} \in \operatorname{span}\{ \psi_{T,k} \}$ such that
\begin{equation*}
  \label{eq:nonSymDtoN}
  \sprod {\psi_{T,k}}{V_T t_{T\ell}} = \sprod  {\psi_{T,k}} {(\tfrac12 I + K_T) u_{T\ell}}
  \qquad
  \forall k = 1, \ldots, |\mathfrak T_\ell(\mathcal F_T)|.
\end{equation*}

This allows us to define the approximate bilinear form
\begin{equation}\label{eq:ApproxBilinearForm}
 \sprod{S_Tu_{T\ell}}{v_{T\ell}} 
 \approx
 \sprod{\widetilde S_Tu_{T\ell}}{v_{T\ell}} 
 :=
 \sprod{t_{T\ell}}{v_{T\ell}},
 \quad \forall u_{T\ell},v_{T\ell} \in \mathfrak W_\ell(\mathcal F_T).
\end{equation}

Let $\underline V_T$, $\underline K_T$ and $\underline{\widetilde S}_T$ denote the matrices
which represent the analytic and approximate bilinear forms induced by the corresponding boundary integral operators 
with respect to the bases $\{\phi_{T,i}\}$ and $\{\psi_{T,k}\}$, i.e.,
\[
  [\underline V_T]_{kl} = \sprod {\psi_{T,k}} {V_T \psi_{T,l}},
  \quad
  [\underline K_T]_{ki} = \sprod {\psi_{T,k}} {K_T \phi_{T,i}},
  \quad
  [\underline{\widetilde S}_T]_{ij} = \sprod {\widetilde S_T\phi_{T,i}} {\phi_{T,j}},
\]
where $i,j$ enumerate the vertices of $\mathfrak T_\ell(\mathcal F_T)$,
and $k, l$ enumerate the triangles of $\mathfrak T_\ell(\mathcal F_T)$.
Furthermore, let $\underline M_T$ denote the mass matrix
\[
  [\underline M_T]_{ki} = \sprod {\psi_{T,k}}{\phi_{T,i}}.
\]
Then, the approximate bilinear form $\sprod{\widetilde S_T \cdot}{\cdot}$ is realized
on our discrete spaces over the element boundary $\partial T$ by the matrix
\[
    \underline{\widetilde S}_T =  \underline M_T^\top\underline V_T^{-1} (\tfrac12 \underline M_T + \underline K_T)
    \in \mathbb R^{\dim \mathfrak W_\ell(\mathcal F_T) \times \dim \mathfrak W_\ell(\mathcal F_T)}.
\]

\subsection{Fully discretized variational problem}
\label{subsec:FullyDiscretizedVariationalProblem}

To obtain the fully discretized variational formulation,
we replace the bilinear form in~\eqref{eq:approxdiscrvf}
by its approximation \eqref{eq:ApproxBilinearForm},
i.e., we seek $u_\ell \in W_\ell$ such that $u_\ell|_{\partial\Omega} = g$ and
\begin{equation}
  \sum_{T \in \mathcal T} \sprod {\widetilde S_T u_{\ell,\partial T}} {v_{\ell,\partial T}} = 0
  \qquad \forall \, v_\ell \in W_{\ell,0}.
  \label{eq:fulldiscrvf}
\end{equation}

For the sake of clarity, we describe in the following how the linear system resulting
from this discretization is constructed.
For an approximate basis function $\varphi_\ell^i$, consider its restriction
$\varphi_{\ell,\partial T}^i \in \mathfrak W_\ell(\mathcal F_T)$ to the boundary of the element $T$.
We denote the coefficients of this restriction with respect to the local nodal piecewise
linear basis functions on $\mathfrak T_\ell(\mathcal F_T)$ by 
$\underline\varphi_{\ell,\partial T}^i \in \mathbb R^{\dim \mathfrak W_\ell(\mathcal F_T)}$.
(If $z_i \notin \partial T$, we have $\underline\varphi_{\ell,\partial T}^i = 0$.)
We gather all these column vectors into a matrix
$\underline D_T \in \mathbb R^{\dim \mathfrak W_\ell(\mathcal F_T) \times |\mathcal N|}$.
Then, the global stiffness matrix is assembled element-wise as
\[
    \underline K = \sum_{T \in \mathcal T} \underline D_T^\top \underline{\widetilde S}_T \underline D_T
    \in \mathbb R^{|\mathcal N| \times |\mathcal N|}.
\]
After eliminating the Dirichlet boundary conditions by means of homogenization in the usual way,
we obtain the final linear system.
At this point we should emphasize that the local auxiliary triangulations $\mathfrak T_\ell(\mathcal F_T)$
are used only to compute the element stiffness matrices. The level of refinement $\ell$
chosen for them has no influence on the size of the global stiffness matrix.

\section{Implementation and numerical examples}
\label{sec:NumericalExperiments}

In this section, we give some implementation details as well as numerical experiments. The computations are done on tetrahedral and polyhedral meshes. For the sake of simplicity, we restrict ourselves to the case of scalar valued diffusion, i.e., $A=\alpha I$ for some $\alpha>0$, and a vanishing reaction term $c=0$. Furthermore, the experiments are carried out with constant and continuously varying convection~$b$. The method is studied for the case of decreasing diffusion $\alpha\to 0$. Standard numerical schemes like the finite element method become unstable when applied to this type of convection-dominated problems. Typically, the issue manifests itself in the form of spurious oscillations.  The critical quantity here is the mesh P\'eclet number
\[
 \mathrm{Pe}_T = \frac{h_T|b_T|}{\alpha_T},\quad T\in\mathcal T,
\]
which should be bounded by $2$ for standard finite element methods. When decreasing the diffusion, the mesh P\'eclet number increases and we expect oscillations. This is due to the fact that the boundary value problem gets closer to a transport equation and thus, boundary layers appear near the outflow boundary.

In addition to stability, we study the number of GMRES iterations which are used to compute the approximate solution of the resulting system of linear equations.

\subsection{Implementation details}
\label{sec:Implementation}

\subsubsection{Preprocessing}
All computations regarding the convection-adapted trial functions can be done in a preprocessing step. In the case of non-constant convection, diffusion and reaction, these terms are first projected into the space of piecewise constant functions over the edges, faces and elements of the mesh. Afterwards, the Dirichlet traces of the trial functions are computed on the edges and faces. Here, an analytic formula is utilized on each edge $E\in\mathcal E$, and subsequently, the two-dimensional convection-diffusion-reaction problems are treated separately on each face $F\in\mathcal F$ according to the SUPG formulation~\eqref{eq:VFonFace}. The stabilization parameter~$\delta_F$ is chosen to be piecewise constant over the auxiliary triangulation~$\mathfrak T_\ell(F)$ on each face $F\in\mathcal F$ with
\[
  \delta_{F,k} = \begin{cases}
                   h_k/2 & \mbox{for } \mathrm{Pe}_{F,k}>2,\\
                   0     & \mbox{else,}
                 \end{cases}
\]
where $k$ enumerates the triangles of~$\mathfrak T_\ell(F)$ and the local P\'eclet number is defined as
\[
 \mathrm{Pe}_{F,k} = \frac{h_k|b_F|}{\alpha_F}.
\]
The auxiliary triangulations~$\mathfrak T_\ell(F)$ of level~$\ell\in\mathbb N_0$ are constructed as described in Section~\ref{sec:Discretization} and visualized in Figure~\ref{fig:AuxiliaryTriangulation}. But, in case of convection-dominated problems on the faces, we decided to move the midpoint of the mesh, created in $\mathfrak T_0(F)$, into the direction of the projected convection vector. Consequently, the auxiliary meshes get adapted to the local problems. This adaptation is inspired by Shishkin-meshes 
\cite{Shishkin:1988},
see also~\cite{KoptevaORiordan2010,Linss:2010,RoosStynesTobiska:2008},
which are graded in such a way that boundary layers are resolved. The solutions of the resulting systems of linear equations, coming from the SUPG formulation, with non-symmetric, sparse matrices are approximated using the GMRES method, see~\cite{SaadSchultz1986}. As the stopping criterion, we use the reduction of the norm of the initial residual by a factor of $10^{-10}$.

Another preprocessing step is the computation of the matrices arising from the local boundary integral formulations. Here, we use the BEM code developed in the PhD thesis by C.~Hofreither~\cite{HofreitherPHD}, which is based on a fully numerical integration scheme described in~\cite{SauterSchwab:eng}. The inversions of the local single layer potential matrices~$\underline{V}_T$ are performed with an efficient LAPACK routine.

\subsubsection{Assembly and solution}
The assembling of the global stiffness matrix is performed element-wise as described in Section~\ref{subsec:FullyDiscretizedVariationalProblem}.
The resulting system of linear equations, which is again sparse and non-symmetric, is treated by GMRES. For the global problem, however, we use the reduction of the norm of the initial residual by a factor of $10^{-6}$ as the stopping criterion. In our numerical experiments, the GMRES iterations are carried out without preconditioning in general. However, we also implemented a simple geometric row scaling (GRS) preconditioner, see~\cite{GordonGordon2010}, i.e., a diagonal preconditioner
\[C^{-1}  = \operatorname{diag}(1/\|\underline{K}_j\|_p),\]
where by $\underline{K}_j$ we mean the $j$-th row of the global stiffness matrix, and we choose the vector norm with $p=1$.

\subsubsection{Improvements}
The proposed method is highly parallelizable, especially the preprocessing steps. The two-dimensional convection-diffusion-reaction problems on the faces are independent of each other, and can thus be treated in parallel. Furthermore, the subsequent setup of the boundary integral matrices and of $\underline{\widetilde S}_T$ can be parallelized on an element level as well. Even the computations of the single entries of each boundary integral matrix are independent of each other.

In the implementation we use another observation to reduce the computational complexity. In the case of constant convection, diffusion and reaction terms, the local boundary integral matrices and the problems on the edges and faces are identical for elements which differ by some translation only. Therefore, we build a look-up table in a preprocessing step such that redundant computations are avoided.

\subsection{Numerical Experiments}
\label{subsec:NumericalExperiments}

\subsubsection{Experiment 1}
\begin{figure}[tb]
 \centering
 \includegraphics[width=0.42\textwidth]{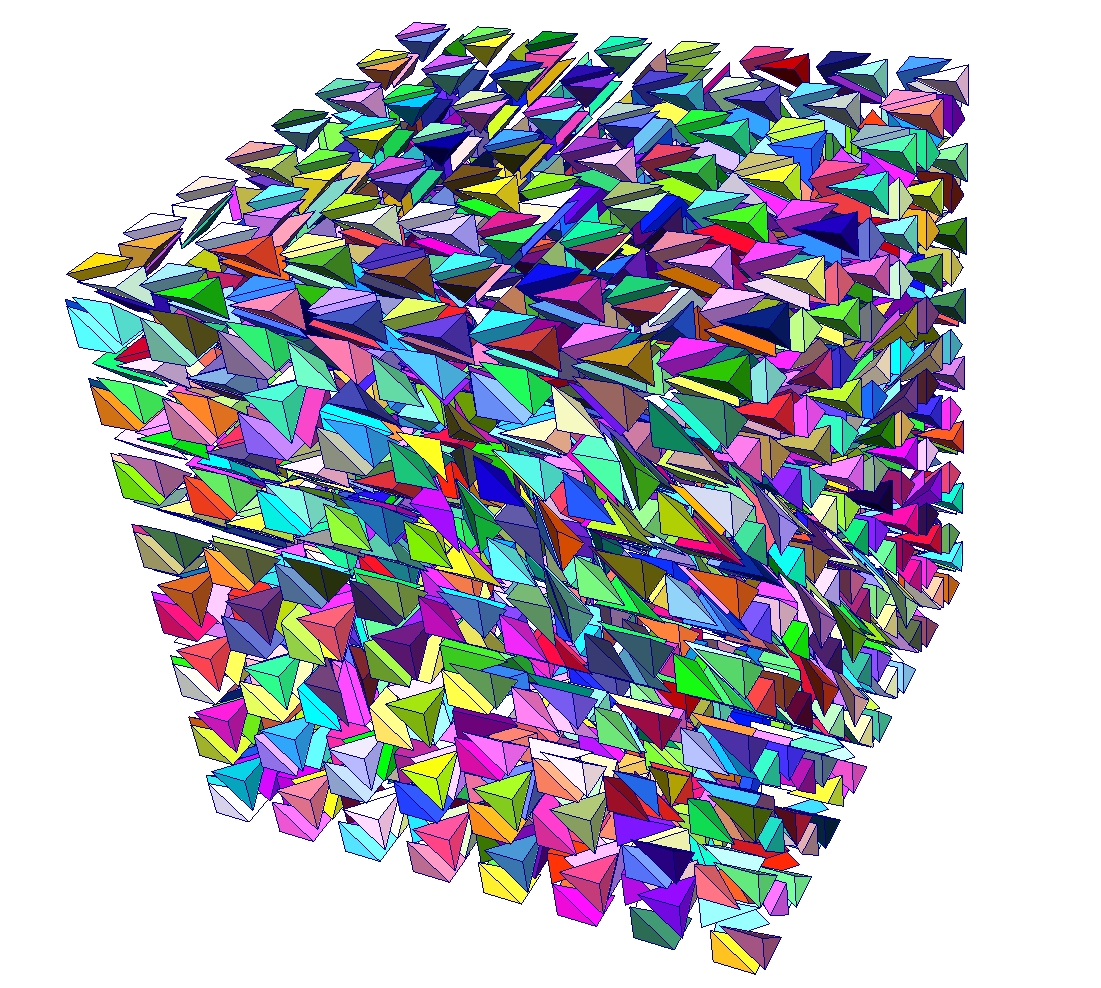}\hfil
 \includegraphics[width=0.42\textwidth]{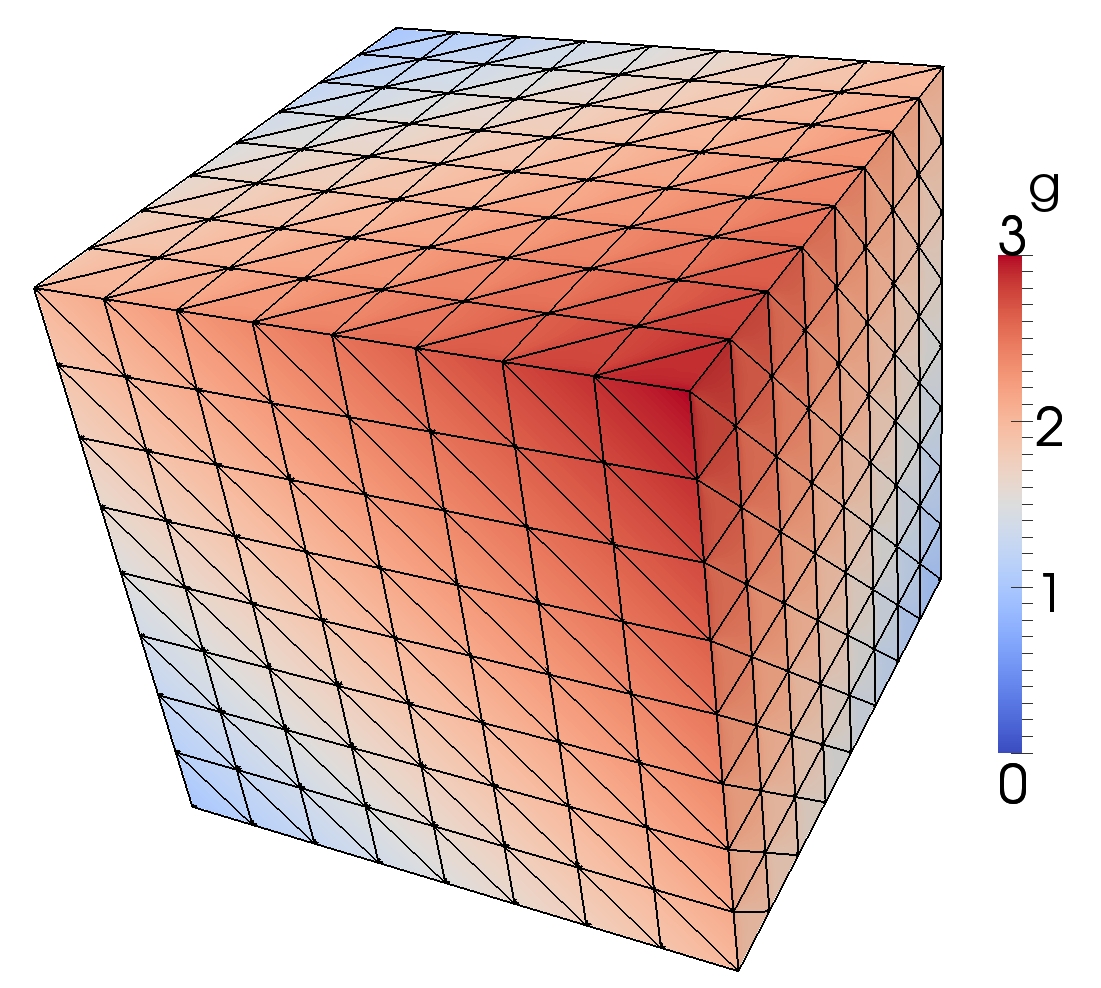}
 \caption{Visualisation of tetrahedral mesh and Dirichlet data for Experiment~1.}
 \label{fig:TetrahedralMesh}
\end{figure}
In the first numerical experiment, a problem with constant convection and diffusion terms is studied. Let $\Omega=(0,1)^3$, and let us consider the boundary value problem 
\begin{equation*}
 \begin{aligned}
  -\alpha\Delta u + b\cdot\nabla u &= 0 &&\mbox{in } \Omega,\\
  u &= g &&\mbox{on } \Gamma,
 \end{aligned}
\end{equation*}
where $b = (1,0,0)^\top$ and $g(x) = x_1+x_2+x_3$.
The domain~$\Omega$ is discretized with tetrahedral elements, see Figure~\ref{fig:TetrahedralMesh}. The mesh consists of 3072 elements, 6528 faces, 4184 edges and 729 nodes of which 343 nodes lie in the interior of~$\Omega$. Thus, the number of degrees of freedom in the BEM-based FEM is equal to 343 in this example. The maximal element diameter is $h_\mathrm{max}\approx 0.22$. The mesh is chosen rather coarse, but it is well suited for the study of stability.

Since the convection and diffusion parameters are constant over the whole domain, the look-up table is applied to speed up the computations. Instead of the before mentioned numbers of geometrical object, we only have to treat 48 elements, 42 faces and 13 edges in the preprocessing step, where the traces of the trial functions are computed and the local stiffness matrices are set up.

To handle the Dirichlet boundary condition, we apply pointwise interpolation of the data~$g$ to obtain an extension onto~$\Gamma_S$. The interpolant is bounded by $0$ from below and by $3$ from above on~$\Gamma$. The convection-diffusion problem satisfies the maximum principle and therefore, we know that $0\leq u\leq 3$ everywhere for the exact solution. To study stability of the BEM-based FEM, the maximum principle is checked for the approximate solution $u_\ell\in W_\ell$ obtained by~\eqref{eq:fulldiscrvf}. Since the trial functions fulfill convection-diffusion problems on the faces and edges and since the maximum principle is also valid there, the maximal values of $u_\ell$ should by reached in the nodes of the mesh. However, because of oscillations coming from the SUPG methods on the faces, the maximal values might be found at some auxiliary node. Consequently, the maximum principle is tested on the whole skeleton~$\Gamma_S$.

\begin{table}[tb]
 \centering
 \caption{Verifying maximum principle in Experiment 1.}
 \label{tab:MaximumPrincipleTet}
 \begin{tabular}{|r|r|rr|rr|rr|}
 \hline & & \multicolumn{2}{c|}{classic FEM} & \multicolumn{4}{c|}{BEM-based FEM} \\
 & & \multicolumn{2}{c|}{linear} & \multicolumn{2}{c|}{straightforward} & \multicolumn{2}{c|}{hierarchical ($\ell=2$)} \\
 $\alpha$ & $\mathrm{Pe}_h$ & $u_\mathrm{min}$ & $u_\mathrm{max}$ & $u_\mathrm{min}$ & $u_\mathrm{max}$ & $u_\mathrm{min}$ & $u_\mathrm{max}$ \\\hline
  $1.0e-1$ & $2$ & $0.00$ & $3.00$ & $0.00$ & $3.00$ & $0.00$ & $3.00$ \\
  $5.0e-2$ & $4$ & $0.00$ & $3.00$ & $0.00$ & $3.00$ & $0.00$ & $3.00$ \\
  $2.5e-2$ & $9$ & $0.00$ & $3.00$ & $0.00$ & $3.00$ & $0.00$ & $3.00$ \\
  $1.0e-2$ & $22$ & $-0.55$ & $3.00$ & $0.00$ & $3.00$ & $-0.01$ & $3.00$ \\
  $5.0e-3$ & $43$ & $-1.14$ & $3.00$ & $0.00$ & $3.00$ & $-0.01$ & $3.00$ \\
  $2.5e-3$ & $87$ & $-1.85$ & $3.07$ & $0.00$ & $3.00$ & $-0.01$ & $3.00$ \\
  $1.0e-3$ & $217$ & $$ & $$ & $0.00$ & $3.00$ & $-0.01$ & $3.00$ \\
  $5.0e-4$ & $433$ & $$ & $$ & $0.00$ & $3.00$ & $-0.01$ & $3.00$ \\
  $2.5e-4$ & $866$ & $$ & $$ & $-142.89$ & $399.06$ & $-0.01$ & $3.00$ \\
  $1.0e-4$ & $2165$ & $$ & $$ & $-68.85$ & $41.00$ & $-0.01$ & $3.00$ \\
  $5.0e-5$ & $4330$ & $$ & $$ & $$ & $$ & $-0.01$ & $3.08$ \\
  $2.5e-5$ & $8660$ & $$ & $$ & $$ & $$ & $-0.01$ & $14.72$ \\\hline
 \end{tabular}
\end{table}
Table~\ref{tab:MaximumPrincipleTet} gives a comparison of the classical Finite Element Method with piecewise linear trial functions and without stabilization, the straightforward BEM-based FEM proposed in~\cite{HofreitherLangerPechstein2011} with linear trial functions on the faces and the new, convection-adapted BEM-based FEM with $\ell=2$. The classical FEM fulfills the discrete maximum principle until $\alpha=2.5e-2$, which corresponds to a P\'eclet number of $9$. The BEM-based strategies, which incorporate the behavior of the differential operator into the approximation space, are more stable. The method in~\cite{HofreitherLangerPechstein2011} passes the test up to $\alpha=5.0e-4$, which corresponds to $\mathrm{Pe}_h=433$. In the new, proposed method we might have oscillations occurring in the approximation of the basis functions fulfilling  convection-dominated problems on the faces. If we neglect these small deviations in the third digit after the decimal point, the proposed method reaches even $\alpha=1.0e-4$, i.e. $\mathrm{Pe}_h=2165$, for $\ell=2$ without violation of the maximum principle.

\begin{table}[tb]
 \centering
 \caption{Verifying maximum principle in Experiment 1 for $\ell=1,2,3$.}
 \label{tab:MaximumPrincipleTetDiffL}
 \begin{tabular}{|r|r|rr|rr|rr|}
 \hline & & \multicolumn{2}{c|}{$\ell=1$} & \multicolumn{2}{c|}{$\ell=2$} & \multicolumn{2}{c|}{$\ell=3$} \\
 $\alpha$ & $\mathrm{Pe}_h$ & $u_\mathrm{min}$ & $u_\mathrm{max}$ & $u_\mathrm{min}$ & $u_\mathrm{max}$ & $u_\mathrm{min}$ & $u_\mathrm{max}$ \\\hline
  $5.0e-3$ & $43$ & $-0.01$ & $3.00$ & $-0.01$ & $3.00$ & $0.00$ & $3.00$ \\
  $2.5e-3$ & $87$ & $-0.01$ & $3.00$ & $-0.01$ & $3.00$ & $0.00$ & $3.00$ \\
  $1.0e-3$ & $217$ & $-0.01$ & $3.00$ & $-0.01$ & $3.00$ & $0.00$ & $3.00$ \\
  $5.0e-4$ & $433$ & $-0.01$ & $3.00$ & $-0.01$ & $3.00$ & $0.00$ & $3.00$ \\
  $2.5e-4$ & $866$ & $-0.01$ & $3.00$ & $-0.01$ & $3.00$ & $0.00$ & $3.00$ \\
  $1.0e-4$ & $2165$ & $-0.01$ & $5.19$ & $-0.01$ & $3.00$ & $0.00$ & $3.00$ \\
  $5.0e-5$ & $4330$ & $-6.72$ & $169.93$ & $-0.01$ & $3.08$ & $0.00$ & $3.00$ \\
  $2.5e-5$ & $8660$ & $-4.8e+6$ & $1.3e+7$ & $-0.01$ & $14.72$ & $0.00$ & $3.00$ \\
  $1.0e-5$ & $21651$ & $$ & $$ & $-7.3e+3$ & $3.1e+4$ & $0.00$ & $26.11$ \\
  $5.0e-6$ & $43301$ & $$ & $$ & $$ & $$ & $-14.89$ & $36.25$ \\\hline
 \end{tabular}
\end{table}
Next, we study the influence of the auxiliary triangulations of the faces on the convection-adapted BEM-based FEM. In Table~\ref{tab:MaximumPrincipleTetDiffL}, the minimal and maximal values $u_\mathrm{min}$ and $u_\mathrm{max}$ of the approximate solution are listed for different levels~$\ell$ of the auxiliary meshes. The higher $\ell$ is chosen, the longer the discrete maximum principle is valid. For $\ell=3$, we even have stability until $\alpha=2.5e-5$, i.e., $\mathrm{Pe}_h=8660$. The enhanced stability can be explained by the improved approximations of the boundary value problems on the edges and faces used to construct the trial functions. Obviously, the local oscillations in the construction of basis functions are reduced such that they have less effect to the global approximation.

\begin{table}[tb]
 \centering
 \caption{Comparing GMRES-iterations in Experiment 1 for the straightforward method and $\ell=1,2,3$.}
 \label{tab:GMRESIterationsTet}
 \begin{tabular}{|r|r|r|r|r|}
 \hline $\alpha$ & s.f. & $\ell=1$ & $\ell=2$ & $\ell=3$ \\\hline
  $ 5.0e-3$ & $ 30$ & $ 28$ & $ 25$ & $ 23$ \\
  $ 2.5e-3$ & $ 33$ & $ 28$ & $ 26$ & $ 24$ \\
  $ 1.0e-3$ & $ 36$ & $ 28$ & $ 26$ & $ 24$ \\
  $ 5.0e-4$ & $ 36$ & $ 28$ & $ 25$ & $ 23$ \\
  $ 2.5e-4$ & $ 311$ & $ 28$ & $ 24$ & $ 23$ \\
  $ 1.0e-4$ & $ 300$ & $ 30$ & $ 24$ & $ 23$ \\
  $ 5.0e-5$ & $ $ & $ 49$ & $ 25$ & $ 23$ \\
  $ 2.5e-5$ & $ $ & $ 302$ & $ 31$ & $ 23$ \\
  $ 1.0e-5$ & $ $ & $ $ & $ 92$ & $ 29$ \\
  $ 5.0e-6$ & $ $ & $ $ & $ $ & $ 55$ \\\hline
 \end{tabular}
\end{table}
In Table~\ref{tab:GMRESIterationsTet}, the numbers of GMRES iterations are given without preconditioning. The GMRES solver for the proposed BEM-based FEM converges faster than for the preceding scheme. For increasing~$\ell$ the convergence slightly improves. Furthermore, the iteration numbers stay bounded without the help of any preconditioning until the maximum principle is violated.

\subsubsection{Experiment 2}
\begin{figure}[tb]
 \centering
 \includegraphics[width=0.42\textwidth]{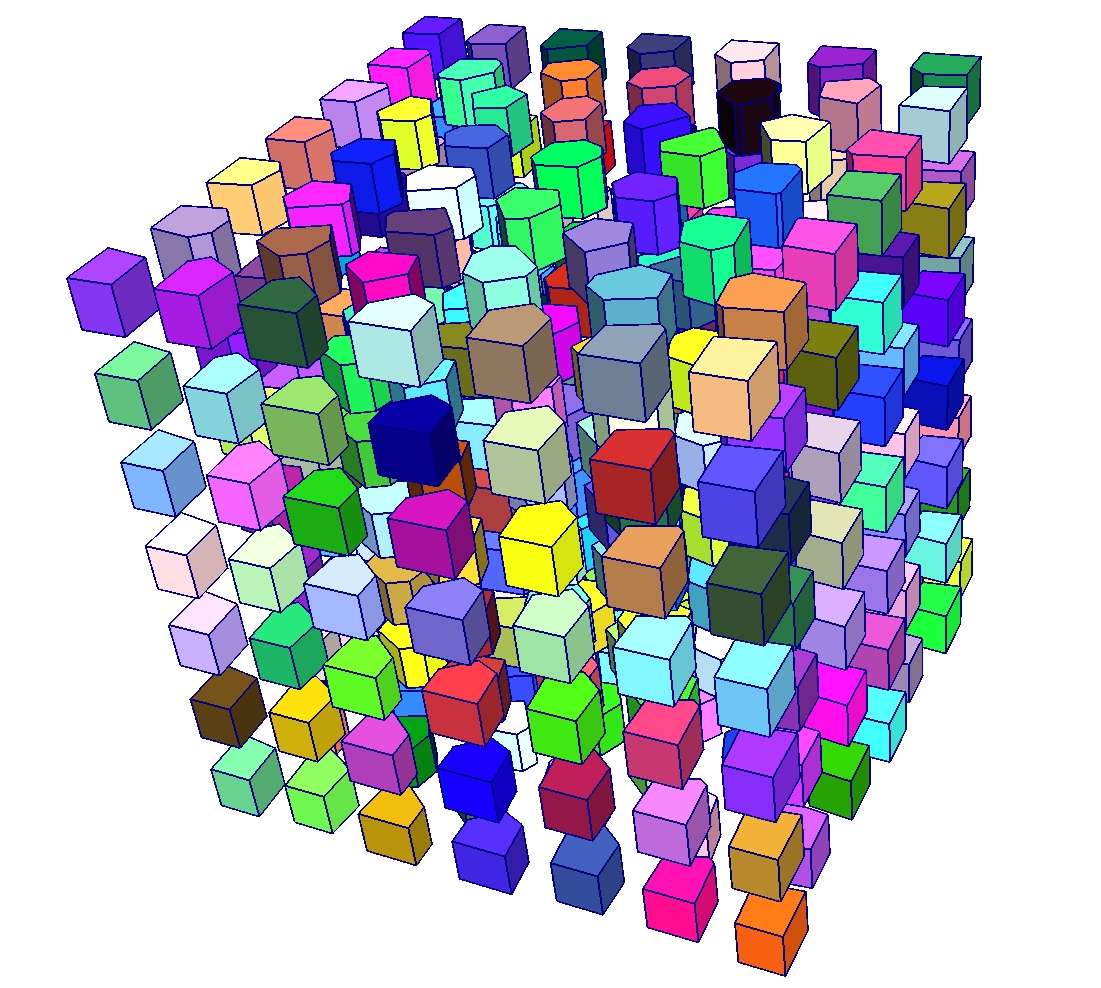}\hfil
 \includegraphics[width=0.42\textwidth]{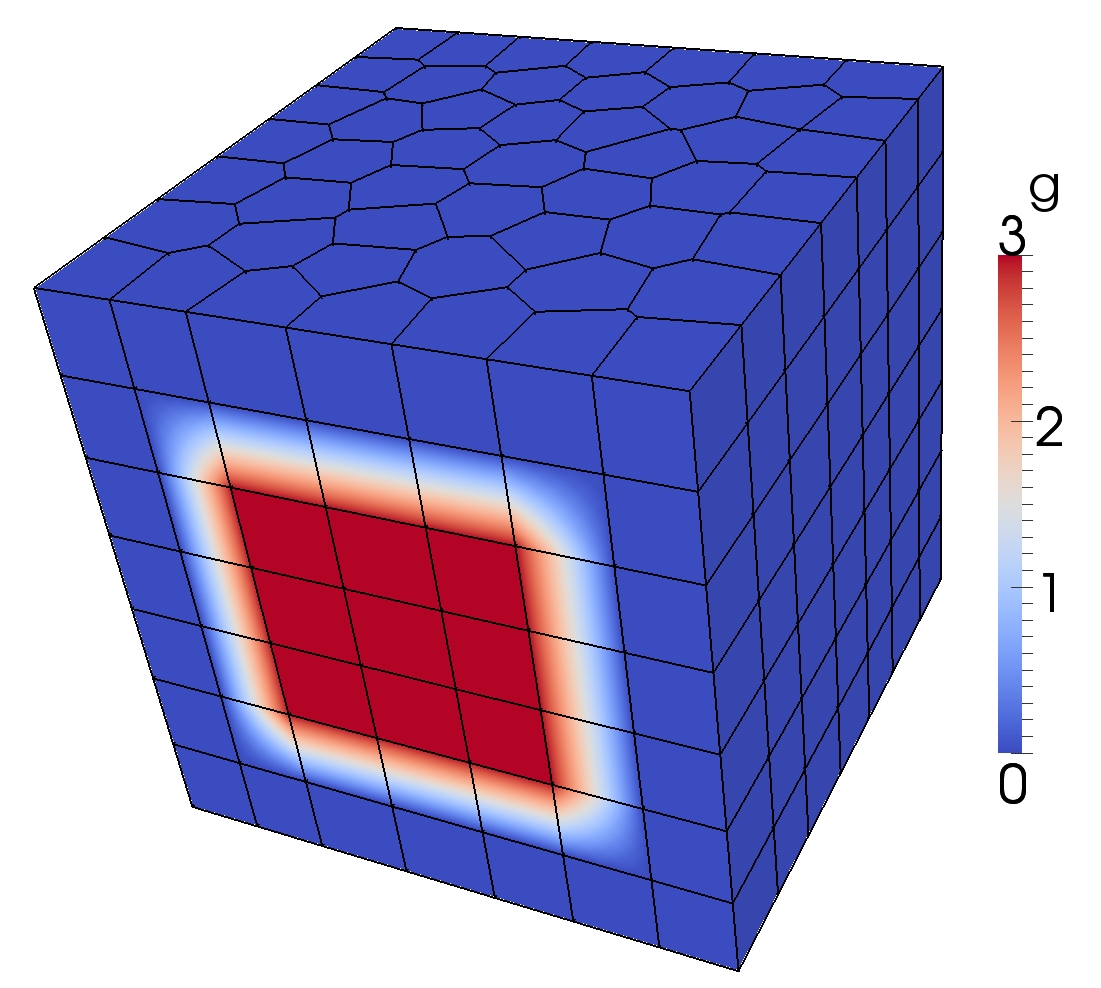}
 \caption{Visualisation of polyhedral mesh and Dirichlet data for Experiment~2.}
 \label{fig:Ex2DirichletData}
\end{figure}
In the next numerical experiment, we consider a convection-diffusion problem with non-constant convection vector. In order to compare the experiments, let $\Omega=(0,1)^3$.
We solve
\begin{equation*}
 \begin{aligned}
  -\alpha\Delta u + b\cdot\nabla u &= 0 &&\mbox{in } \Omega,\\
  u &= g &&\mbox{on } \Gamma,
 \end{aligned}
\end{equation*}
where
\[
 b(x) = \frac{0.85}{\sqrt{(1-x_1)^2+(1-x_3)^2}}\left(\begin{array}{c}x_3-1\\ 0\\ 1-x_1\end{array}\right)
\]
and $g$ is chosen such that it is piecewise bilinear and continuous with $0\leq g\leq 3$ on one side of the unit cube and zero on all others, see Figure~\ref{fig:Ex2DirichletData}.
The convection vector $b$ is scaled in such a way that the P\'eclet numbers in the computations are comparable with those of Experiment~1.
The convection is a rotating field around the upper edge of the unit cube $\Omega$, which lies in the front when looking at Figure~\ref{fig:Ex2DirichletData}. 
Consequently, we expect that the non-zero Dirichlet data is transported towards the upper side of the cube for low diffusion.

This time, the domain~$\Omega$ is decomposed into prisms having general polygonal ends, see Figure~\ref{fig:Ex2DirichletData}. 
The polyhedral mesh consists of 350 elements, 1450 faces, 1907 edges and 808 nodes of which 438 nodes lie in the interior of~$\Omega$. 
Thus, the number of degrees of freedom in the BEM-based FEM is equal to 438. 
The maximal diameter of the elements is $h_\mathrm{max}\approx 0.25$ and the discretization was chosen such that $h_\mathrm{max}$ is approximately the same as in Experiment~1.

In our experiments, the polyhedral mesh has less elements, faces and edges than the tetrahedral discretization. 
This is beneficial concerning the computations in the preprocessing step. 
Less local problems have to be solved on edges and faces and there are less boundary element matrices which have to be set up. 
Furthermore, polyhedral discretizations admit a high flexibility while meshing complex geometries.

\begin{table}[tb]
 \centering
 \caption{Verifying maximum principle in Experiment 2 for $\ell=2$ and number of iterations with/without preconditioning.}
 \label{tab:MaximumPrinciplePolyDiffL}
 \begin{tabular}{|r|r|rr|c|c|}
 \hline$\alpha$ & $\mathrm{Pe}_h$ & $u_\mathrm{min}$ & $u_\mathrm{max}$ & iter. & iter. (prec.) \\\hline
  $1.0e-1$ & $2$ & $0.00$ & $3.00$ & $20$ & $20$ \\
  $5.0e-2$ & $4$ & $0.00$ & $3.00$ & $20$ & $21$ \\
  $2.5e-2$ & $9$ & $0.00$ & $3.04$ & $20$ & $21$ \\
  $1.0e-2$ & $22$ & $0.00$ & $3.07$ & $23$ & $22$ \\
  $5.0e-3$ & $43$ & $-0.01$ & $3.26$ & $29$ & $23$ \\
  $2.5e-3$ & $86$ & $-0.04$ & $3.37$ & $42$ & $24$ \\
  $1.0e-3$ & $216$ & $-0.10$ & $3.38$ & $45$ & $23$ \\
  $5.0e-4$ & $431$ & $-0.13$ & $3.45$ & $48$ & $22$ \\
  $2.5e-4$ & $863$ & $-0.15$ & $3.51$ & $51$ & $21$ \\
  $1.0e-4$ & $2157$ & $-0.15$ & $3.53$ & $52$ & $21$ \\
  $5.0e-5$ & $4313$ & $-0.16$ & $3.57$ & $58$ & $23$ \\
  $2.5e-5$ & $8627$ & $-0.25$ & $4.38$ & $69$ & $28$ \\\hline
 \end{tabular}
\end{table}
In Table~\ref{tab:MaximumPrinciplePolyDiffL}, we list the minimal and maximal values of the approximation $u_\ell$ on the skeleton for $\ell=2$ to verify the discrete maximum principle.
Furthermore, the numbers of GMRES iterations are given with and without preconditioning.

The first observation is that the number of GMRES iterations increases when the diffusion~$\alpha$ tends to zero. Thus, the iteration count is not bounded in this experiment. 
However, this behavior correlates with the violation of the maximum principle and is therefore the result of inaccuracies.
Already with the help of the simple geometric row scaling preconditioner, we overcome the increase of the iteration number.

A more detailed discussion is needed for the discrete maximum principle. 
In Table~\ref{tab:MaximumPrinciplePolyDiffL}, we observe that this principle is violated in a relatively early stage for $\alpha=2.5e-2$, which corresponds to $\mathrm{Pe}_h=9$. 
However, the increase of $u_\mathrm{max}$ and the decrease of $u_\mathrm{min}$ is fairly slow for increasing P\'eclet number. 

Here, one has to point out that the computations are done on a polyhedral mesh with a globally continuous approximation $u_\ell$. 
This, by itself, is a current field of research even without dominant convection, see~\cite{BeiraoBrezziMariniRusso2014preprint}. 
The geometry of polygonal faces is more complex than the triangles in Experiment~1, and thus, the computations on the faces are more involved.

\begin{figure}[tb]
 \centering
 \includegraphics[width=0.42\textwidth]{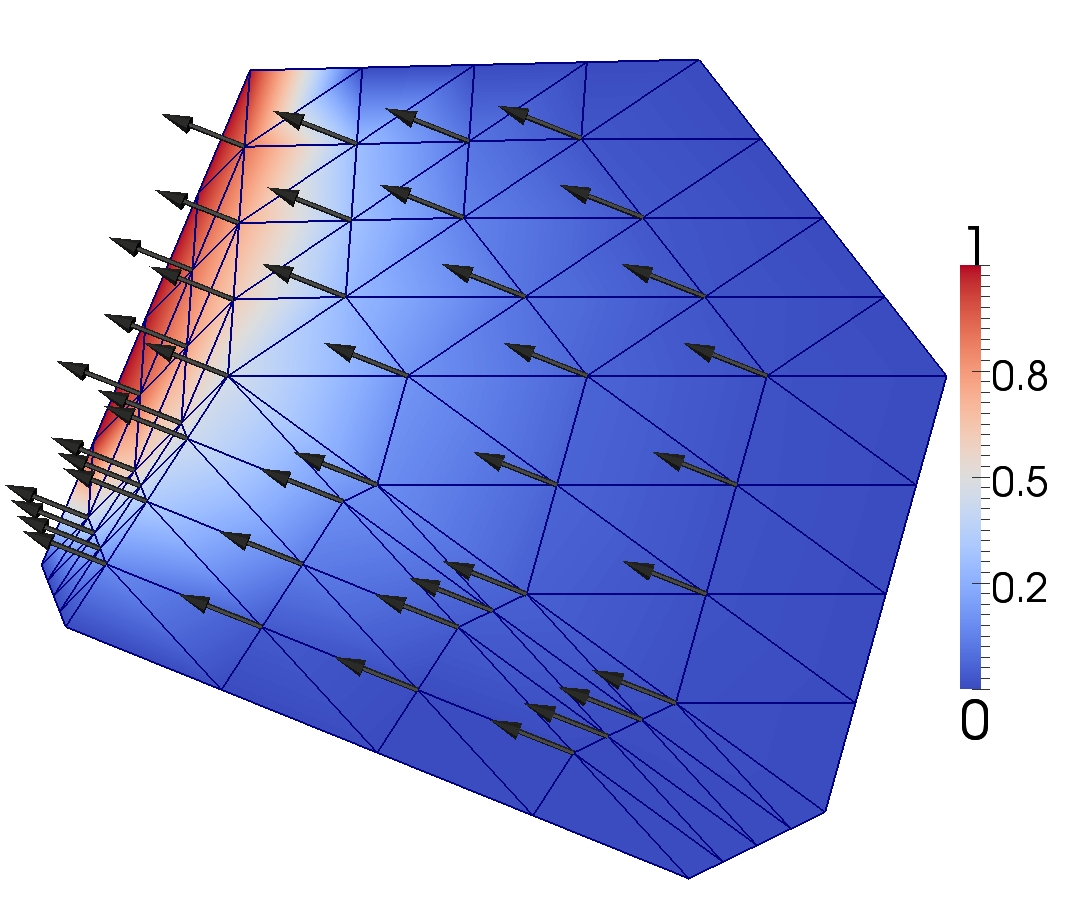}\hfil
 \includegraphics[width=0.42\textwidth]{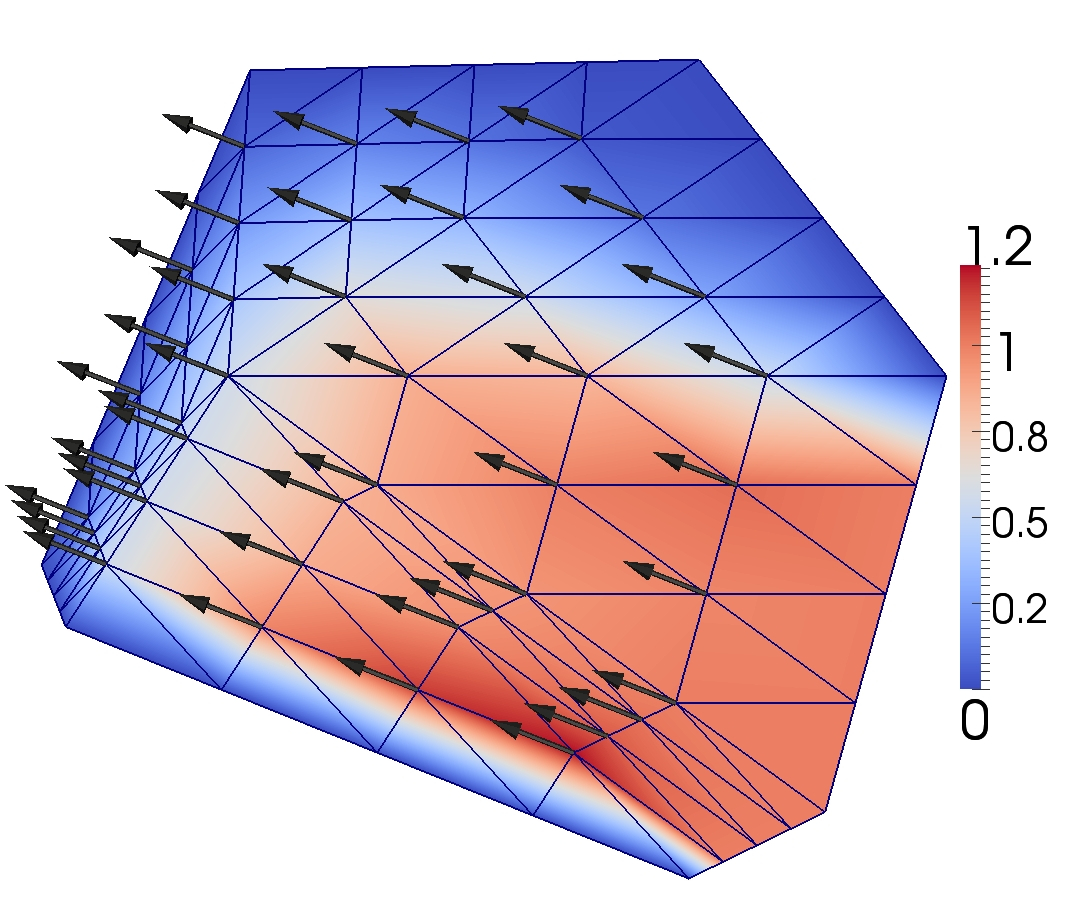}
 \caption{Approximations of basis functions on polygonal face, projected convection vector and auxiliary triangulation with appropriately (left) and not appropriately (right) resolved boundary layer.}
 \label{fig:ApproxBasisFunc}
\end{figure}
Figure~\ref{fig:ApproxBasisFunc} shows the approximation of two different basis functions over the same polygonal face, the auxiliary triangulation and the projected convection vector. 
We can see how the local mesh has been adapted to the underlying differential operator, namely by moving the node, which lay initially in the center of the polygon, into the direction of the convection.
In certain constellations, the boundary layers are not resolved appropriately. 
In the left picture of Figure~\ref{fig:ApproxBasisFunc}, the approximation of the basis function is satisfactory. 
In the right picture, however, oscillations occur in the lower right corner due to the relatively large triangles near the boundary. 
In many cases these situations are already resolved quite well by the simple mesh adaptation. 
When we introduced the moving of the auxiliary nodes in the implementation, the numerical results improved.
Thus, we expect that a better adaptation of the local meshes, and consequently a better approximation of the local problems, improves the stability of the BEM-based FEM such that we would obtain comparable results to Experiment~1 for the discrete maximum principle.

\begin{figure}[tb]
 \centering
 \includegraphics[width=0.42\textwidth]{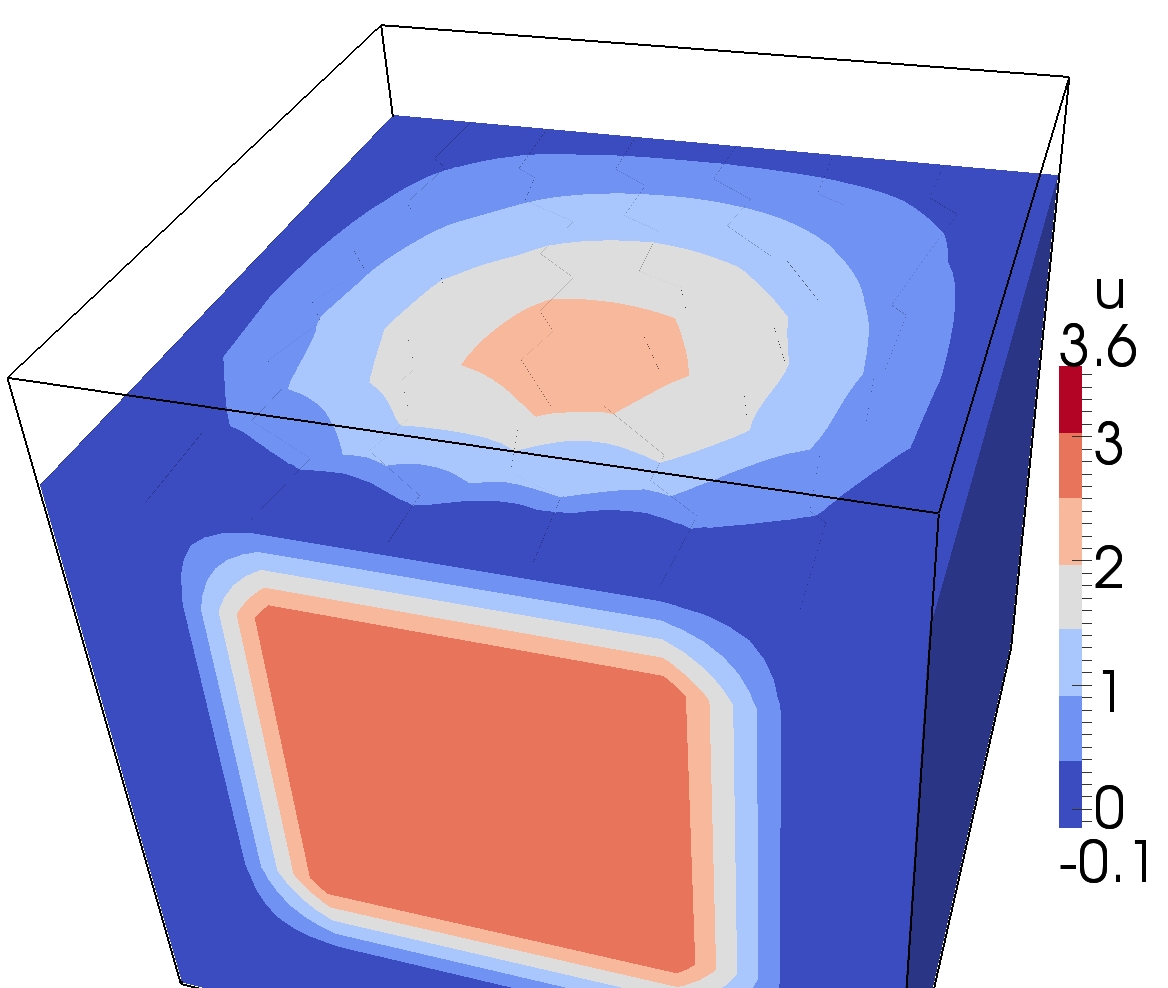}\hfil
 \includegraphics[width=0.42\textwidth]{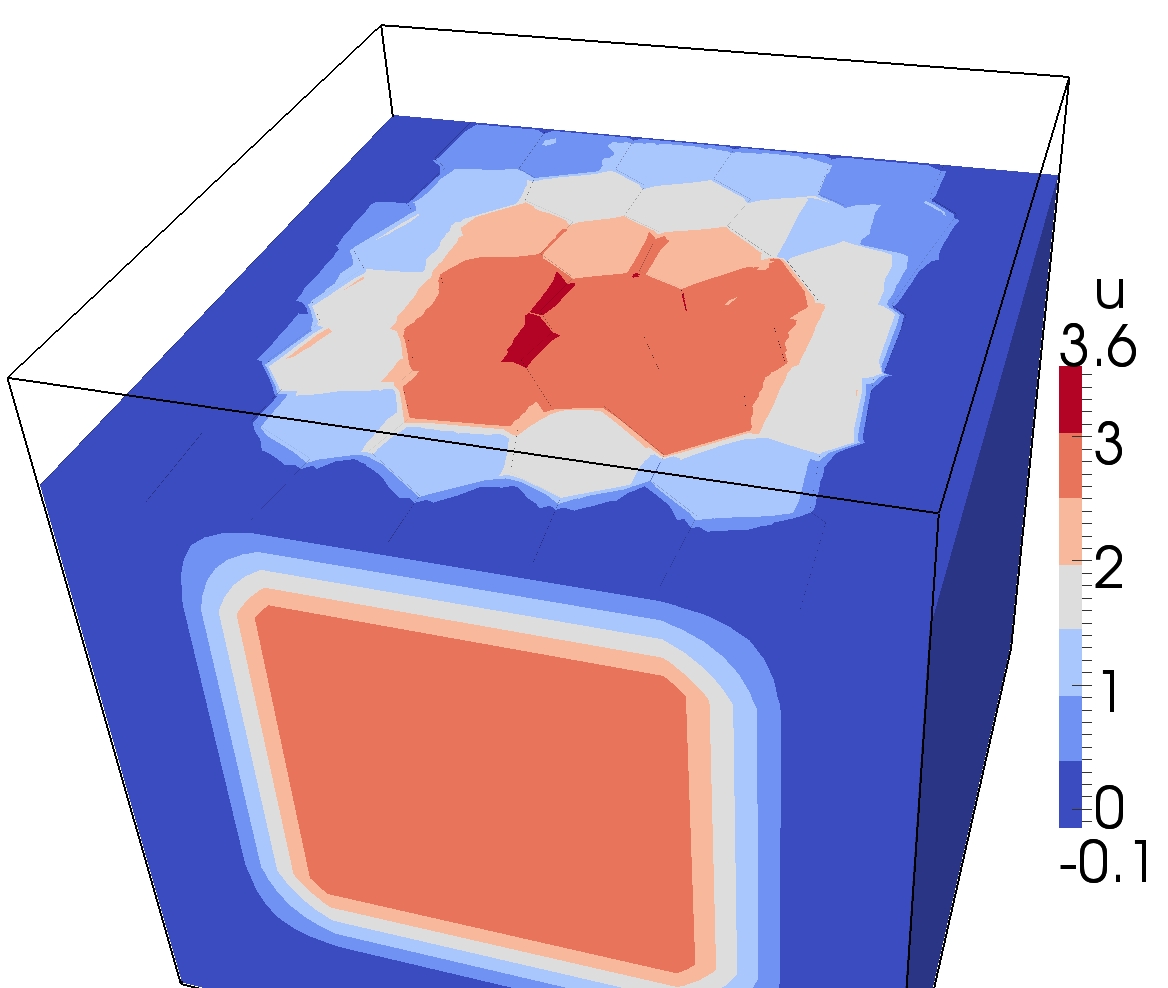}
 \caption{Cut through the domain~$\Omega=(0,3)^3$ and visualisation of the approximation in Experiment 2 for $\alpha=2.5e-2$ (left) and $\alpha=5.0e-5$ (right).}
 \label{fig:Approx}
\end{figure}
Finally, in Figure~\ref{fig:Approx}, the approximation $u_\ell$ is visualized for $\ell=2$ and two different values of diffusion $\alpha=2.5e-2$ and $\alpha=5.0e-5$. 
The domain~$\Omega$ has been cut through, such that the approximation is visible on a set of polygonal faces which lie in the interior of the domain. 
The expected behavior of the solution can be observed. 
The Dirichlet data is transported into the interior of the domain along the convection vector. 
In the case of the convection-dominated problem, oscillations appear near the outflow boundary.

\section{Conclusion}
\label{sec:Conclusion}
We have derived  new convection-adapted BEM-based FEM discretization schemes
for convection-diffusion-reaction boundary value problems 
that considerably extend the range of applicability 
with respect to the strength of convection.
The numerical results have not only confirmed 
this enhanced stability property of the discretization scheme,
but have also indicated faster convergence of the GMRES solver  
in comparison with the original BEM-based FEM scheme
presented in \cite{HofreitherLangerPechstein2011,HofreitherPHD}.


\end{document}